\documentclass{article}

\usepackage{arxiv}

\usepackage[utf8]{inputenc} 
\usepackage[T1]{fontenc}    
\usepackage{hyperref}       
\usepackage{url}            
\usepackage{booktabs}       
\usepackage{amsfonts}       
\usepackage{nicefrac}       
\usepackage{microtype}      
\usepackage{lipsum}
\usepackage{amsmath}

\usepackage{amsfonts,amssymb,stmaryrd}
\usepackage{graphicx}
\usepackage{amsthm}
\usepackage{xcolor}

\newcommand{\R}{\mathbb{R}}

\newtheorem{theorem}{Theorem}[section]
\newtheorem{proposition}{Proposition}[section]

\newtheorem{corollary}{Corollary}[section]

\newtheorem{remark}{Remark}[section]

\usepackage{subcaption}
\usepackage{cite}
\usepackage{authblk}
\newcommand{\p}{\partial}
\newcommand{\bb}{\begin{equation}}
\newcommand{\ee}{\end{equation}}
\newcommand{\ba}{\begin{array}}
\newcommand{\ea}{\end{array}}
\newcommand{\f}{\frac}
\usepackage[all]{xy}
\newcommand{\ds}{\displaystyle}

\newcommand{\al}{\alpha}

\newcommand{\I}{{\cal I}_0}

\numberwithin{equation}{section}

\usepackage{ulem}

\title{Persistence properties of a Camassa-Holm type equation with $(n+1)-$order non-linearities}

\author{
Igor~Leite~Freire\thanks{igor.freire@ufabc.edu.br or igor.leite.freire@gmail.com}\\
$^1$Centro de Matem\'atica, Computa\c{c}\~ao e Cogni\c{c}\~ao\\
Universidade Federal do ABC\\
Santo Andr\'e, Brazil}
 
\begin{document}
\maketitle

\begin{abstract}
Lower order conservation laws and symmetries of a family of hyperbolic equations having the Camassa-Holm equation as a particular member are obtained. We show that the equation has two conservation laws with zeroth order characteristics and that its symmetries are generated by translations in the independent variables and a certain scaling, as well as some invariant solutions are studied. Next, we consider persistence and asymptotic properties for the solutions of the equation considered. In particular, we analyse the behaviour of the solutions of the equation for large values of the spatial variable. We show that if the initial data has a certain asymptotic exponential decaying, then such property persists for any time as long as the solution exists. Moreover, depending on the behaviour of the initial data for large values of the spatial variable and if for some further time the solution shares the same behaviour, then it necessarily vanishes identically. Finally, we prove unique continuation results for the solutions of the equation.
\end{abstract}

{\bf MSC classification 2010:} 35A01, 74G25, 37K40, 35Q51.


\newpage
\section{Introduction}\label{sec1}

In \cite{hak-cpde,hak-ijns} Hakkaev and Kirchev considered the following equation
\bb\label{1.0.1}
u_t-u_{txx}+\f{(n+1)(n+2)}{2}u^nu_x=\f{n(n-1)}{2}u^{n-2}u_x^3+2n u^{n-1}u_xu_{xx}+u^n u_{xxx},
\ee
where $n$ is a constant, that henceforth will be assumed to be a positive integer, and $u=u(t,x)$. We observe that whenever $n=1$ we recover the celebrated Camassa-Holm (CH) equation \cite{chprl}, which is nowadays a by far famous equation not only for its relevance in hydrodynamics \cite{chprl,const-arc}, but also for its rich mathematical properties, see \cite{chprl,const2000-1,const-jmp,const-arc,raspajde,raspaspam,freire-jpa,freire-arxiv,henry-jnmp,henri-na,henry-jhde} and references therein about several properties of the CH equation. For the CH equation the variables $t$ and $x$ are regarded as time and space, respectively, so that we maintain this terminology for \eqref{1.0.1}.

In regard to \eqref{1.0.1}, it is locally well-posed in Sobolev spaces $H^s(\R)$, with $s>3/2$, see \cite[Theorem 2]{lai-jde} and \cite[Theorem 3.4]{hak-cpde} (see also \cite[Theorem 2.2]{hak-ijns}), and for $s\in(1,3/2]$ the existence of weak solutions is granted by \cite[Theorem 1]{lai-jde}.

In addition to the results in \cite{hak-cpde,lai-jde}, in \cite{mi-jmaa} was shown the local well-posedness of solutions of \eqref{1.0.1} in Besov spaces (see Theorem 1.1 in the mentioned reference) and under certain circumstances, the solutions are also analytic on a certain open set of $\R^2$, see \cite[Theorem 1.2]{mi-jmaa}.

In \cite{zhou-jmaa} some persistence properties for \eqref{1.0.1} were considered. In particular, among its results there is one concerned with asymptotic properties of \eqref{1.0.1}: in case $n$ is odd and for each $t$ fixed, then $u(t,x)\sim u_0(x)+c(x)e^{-x}$, where $c(x)\rightarrow0$ as $x\rightarrow\infty$, see \cite[Theorem 1.3]{zhou-jmaa}.

The purpose of this paper is to shed light on:
\begin{itemize}
    \item Invariance properties of \eqref{1.0.1}, namely, local diffeomorphisms preserving solutions of the equation, as well as conserved currents and their corresponding conserved quantities;
    \item Unique continuation results for \eqref{1.0.1};
\end{itemize}

In addition to these points, we revisit some results proved in \cite{zhou-jmaa} from a different perspective influenced by recent developments considered in \cite{linares,freire-jpa,raspa-mo,freire-arxiv} and extending to \eqref{1.0.1} some results proved in the mentioned references. To do this, we consider families of functions obtained from a known solution $u(t,x)$ of \eqref{1.0.1}. More precisely, we consider the families
$${\cal F}_1=(f_t)_{t\in[0,T)}\quad\text{and}\quad {\cal F}_2=(F_t)_{t\in[0,T)},$$
where
$$
f_t(x):=\f{n}{2}u(t,x)^{n-1}u_x(t,x)^2+\f{n(n+3)}{2(n+1)}u(t,x)^{n+1}\quad\text{and}\quad F_t(x):=\p_x\Lambda^{-2}\ast f_t(x).
$$

The meaning of $\p_x\Lambda^{-2}$ will be given in Section \ref{sec2}, where we present the notations and conventions used throughout the paper, as well as our main results and outline of the work.

\section{Main results of the paper, their preliminary discussion, and its outline}\label{sec2}

We introduce and fix the notation, notions and conventions used throughout the paper. Next, we state our main results and the manuscript's outline.

\subsection{Notation and conventions} 

By $\|\cdot\|_{p}$ and $\|\cdot\|_{H^s(\R)}$ we denote the usual norms of the Banach spaces $L^p(\R)$, $1\leq p\leq\infty$, and the Sobolev space $H^s(\R)$, $s\in\R$, respectively. The convolution of two functions $f$ and $g$ is denoted by $f\ast g$. If $u$ is a function of two variables $t$ and $x$, that is, $u=u(t,x)$, the function $x\mapsto u(0,x)$ will be frequently denoted by $u_0(x)$, whereas $u_t$ and $u_x$ denote the derivatives of $u$ with respect to the first and the second arguments, respectively. The Helmholtz operator $1-\p_x^2$ and its inverse $(1-\p_x^2)^{-1}$ are denoted by $\Lambda^2$ and $\Lambda^{-2}$, respectively. In particular, $\Lambda^{-2}(f)=g\ast f$, where $g$ is the Green function 
\bb\label{2.1.1}
g(x)=\f{e^{-|x|}}{2}.
\ee

We say that $|f(x)|\sim O(e^{a x})$ as $x\nnearrow\infty$ if there exists some $L$ such that
$$
\lim_{x\rightarrow\infty}\f{|f(x)|}{e^{a x}}=L,
$$
whereas $|f(x)|\sim o(e^{a x})$ as $x\nnearrow\infty$ if
$$
\lim_{x\rightarrow\infty}\f{|f(x)|}{e^{a x}}=0.
$$

The first definition does not exclude the case $L=0$. If this happens, both definitions coincide, meaning that if $|f(x)|\sim o(e^{a x})$ as $x\nnearrow\infty$, then $|f(x)|\sim O(e^{a x})$ as $x\nnearrow\infty$, but the converse may not be true.

Finally, we recall the Gr\"onwall inequality, e.g., see \cite[page 29]{taylor}, which says that if $M(t)$ is non-negative and $t\geq0$ and
$$u(t)\leq A+\int_{0}^tM(s)u(s)ds,$$
then $u(t)\leq A e^{\int_0^tM(s)ds}$. Note that sometimes the original inequality is replaced by $u'(t)\leq M(t)u(t)$.

We recall that equation \eqref{1.0.1} can be rewritten as 
\bb\label{2.1.2}
u_t+u^nu_x+\p_x(1-\p_x^2)^{-1}\left(\f{n}{2}u^{n-1}u_x^2+\f{n(n+3)}{2(n+1)}u^{n+1}\right)=0.
\ee

The results proved in \cite{hak-cpde,mi-jmaa,lai-jde} show that if if $u=u(t,x)$ is a solution of \eqref{1.0.1} such that $u(0,x)=u_0(x)$, for some $u_0\in H^s(\R)$, with $s>3/2$, then the corresponding Cauchy problem has a unique solution $u\in C^0([0,T);H^s(\R))\cap C^1([0,T);H^{s-1}(\R))$, for some $T>0$. Unless otherwise mentioned, throughout this paper by a solution $u$ of \eqref{1.0.1} we mean a function $u\in C^0([0,T);H^s(\R))$, $s>3/2$ (the exception is in some parts of Subsection \ref{subsec2.2}). This fact will be widely used henceforth without further mention.

\subsection{Main results and outline of the paper}\label{subsec2.2}

Our first result is concerned with conserved currents, which are vector fields $(C^0,C^1)$ whose divergence $\p_t C^0+\p_xC^1$ vanishes on the solutions of the equation (for further details, see \cite[Sec. 4.3]{olverbook}). A characteristic of order $m$ ($m=0,\,1,\,2,\cdots$) of a conservation law of an equation $E=0$ is a function $Q$ of $t,x$, $u$ and derivatives of $u$ up to order $m$ such that $\p_t C^0+\p_x C^1=QE$. Therefore, as long as $u$ is a solution of $E=0$, we then have $\p_t C^0+\p_x C^1=0$.

\begin{theorem}\label{teo2.1}
Equation \eqref{1.0.1} admits only two conserved currents with zeroth characteristic $Q_1=1$ and $Q_2=u$, for any positive integer $n$.
\end{theorem}

\begin{corollary}\label{cor2.1}
Assume that $u$ is a solution of \eqref{1.0.1} vanishing as $|x|\rightarrow\infty$, whose derivatives up to second order are bounded. Then the quantities
\bb\label{2.2.1}
{\cal H}_1(t):=\int_\R u(t,x)dx
\ee
and 
\bb\label{2.2.2}
{\cal H}(t):=\int_\R\f{u(t,x)^2+u_x(t,x)^2}{2}dx=\f{1}{2}\|u(t,\cdot)\|_{H^1(\R)}^2
\ee
are constant. In particular, if $u$ is either non-negative or non-positive, then $\|u(t,\cdot)\|_{1}$ is constant. 
\end{corollary}

Note that solutions belonging to $C^0([0,T);H^s(\R))$, with $s>3/2$, will have both \eqref{2.2.1} and \eqref{2.2.2} as conserved quantities, meaning that ${\cal H}_1(t)={\cal H}_1(0)$ and ${\cal H}(t)={\cal H}(0)$, for any $t\in(0,T)$.

Our next result gives the Lie point symmetries of \eqref{1.0.1}. We recall that a Lie point symmetry of a given differential equation is a local diffeomorphism that maps a (strong or smooth) solution of the equation into another solution in a smooth way. As such, the symmetries are fluxes of certain linear operators, their generators, and therefore, once we know the generators we know the symmetries and vice-versa. Geometrically, or even physically, the Lie symmetries can be seen as smooth motions or deformations of the solutions of the equation. 

\begin{theorem}\label{teo2.2}
The Lie point symmetries of the equation \eqref{1.0.1} are
\bb\label{2.2.3}
{\bf v}_1=\p_t,\,\,{\bf v}_2=\p_x,\,\,{\bf v}_3=n t\p_t-u\p_u.
\ee
\end{theorem}

The proof of Theorem \ref{teo2.2} is straightforward, but long. It consists in applying the algorithm for finding the coefficients of the generator of symmetries \cite[Chap. 2]{olverbook}. This will lead us to obtain an over-determined system of linear differential equations to the unknown coefficients whose solution, substituted into the generator of symmetries, will correspond to a linear combination of \eqref{2.2.3}. Since this process is very long, and can be carried out by several computational packages available, e.g., see \cite{stelios1,stelios2}, we omit its demonstration.

The fluxes of the generators \eqref{2.2.3} correspond to translations in time $(t,x,u)\mapsto(t+\epsilon,x,u)$, space $(t,x,u)\mapsto(t,x+\epsilon,u)$ and the scaling $(t,x,u)\mapsto(e^{n \epsilon}t,x,e^{-\epsilon}u)$, where $\epsilon\in\R$. The invariance under the scaling leads us to solutions of the type
$$
u(t,x)=\f{f(x)}{t^{1/n}}.
$$

Substituting this $u$ into \eqref{1.0.1} we obtain the following ODE to $f$:
$$
-\f{1}{n}(f-f'')+\f{(n+1)(n+2)}{2}f^n f'=\f{n(n-1)}{2}f^{n-2}(f')^3+2n f^{n-1}f'f''+f^n f''',
$$
where we omitted the dependence with respect to $x$ and the prime $'$ denotes derivative with respect the independent variable. Two straightforward solutions to this ODE are $f_+(x)=e^x$ and $f_-(x)=e^{-x}$ and, therefore, it is easy to check that they provide solutions of \eqref{1.0.1} that does not conserve \eqref{2.2.1} nor \eqref{2.2.2}. On the other hand, the invariance under translations in $t$ and $x$ naturally lead us to the investigation of travelling waves admitted by \eqref{1.0.1}. 

Let $z=x-ct$ and $u(t,x)=\phi(z)$. Substituting this $u$ into \eqref{1.0.1} and rearranging the result, we have:
\bb\label{2.2.4}
\f{d}{dz}\left(-c\phi+c\phi''+\f{n+2}{2}\phi^{n+2}\right)=\f{d}{dz}\left(\f{n}{2}\phi^{n-1}(\phi')^2+\phi^n\phi''\right),
\ee
where the dependence is again omitted and the prime denotes usual derivative. Integrating the equation above and denoting the constant of integration by $k_1/2$, we have
$$
-c\phi +c\phi''+\f{n+2}{2}\phi^{n+1}=\f{n}{2}\phi^{n-2}(\phi')^2+\phi^n\phi''+\f{k_1}{2}.
$$
Multiplying the latter equation by $\phi'$, integrating once more, denoting the new constant of integration by $k_2/2$ and rearranging the terms, we obtain
\bb\label{2.2.5}
(\phi')^2=\f{c\phi^2-\phi^{n+2}+k_1\phi+k_2}{c-\phi^n}.
\ee

Equation \eqref{2.2.5} shows that, in principle, we should expect problems at the derivatives of $\phi$ if there exists a point $z_0$ such that $\phi(z)\rightarrow c^{1/n}$ as $z\rightarrow z_0$. However, if such a pole is removable, we can hope for the emergence of a solution with some issues regarding differentiability, meaning that we may expect that \eqref{1.0.1} has solutions in the distributional sense. In particular, if the constants of integration $k_1$ and $k_2$ are taken as $0$, then $(\phi')^2=\phi^2$, and as a consequence we have the solution $\phi(z)=c^{1/n}e^{-|z|}$, implying that
\bb\label{2.2.6}
u(t,x)=c^{1/n}e^{-|x-ct|}
\ee
solves \eqref{1.0.1} in the weak sense. Such a solution, called ${\it peakon}$, is orbitally stable, see \cite{hak-ijns}, and is one of the features of the CH equation shared by \eqref{1.0.1}.

In \cite{him-cmp} was shown that if an initial data of the CH equation, for large values of $x$, behaves like the peakon solution, then this property persists for the corresponding solution. Therefore, it is natural to investigate whether the solutions of \eqref{1.0.1} would also inherit the same property using the approach in \cite{him-cmp}.

Our next two results shed light to the above question. In what follows, $T_0\in(0,T)$ is an arbitrary value for which the solutions is defined and we denote by $\I$ the compact set $[0,T_0]$.

\begin{theorem}\label{teo2.3}
Assume that $u\in C^0(\I;H^s(\R))$, $s>3/2$, is a solution of \eqref{1.0.1} with initial condition $u(0,x)=u_0(x)$. Assume that for some $\theta\in(0,1)$, 
    \bb\label{2.2.7}
    |u_0(x)|\sim O(e^{-\theta x})\quad\text{and}\quad |u_0'(x)|\sim O(e^{-\theta x})\quad\text{as}\quad x\nnearrow\infty.
    \ee
Then 
    \bb\label{2.2.8}
    |u(t,x)|\sim O(e^{-\theta x})\quad\text{and}\quad |u_x(t,x)|\sim O(e^{-\theta x})\quad\text{as}\quad x\nnearrow\infty,
    \ee
uniformly in $\I$.
\end{theorem}

Theorem \ref{teo2.3} could be inferred from \cite[Theorem 1.2]{zhou-jmaa}. However, the approach used in \cite{zhou-jmaa} is different of ours. Moreover, Theorem \ref{teo2.3} is a {\it sine qua non} ingredient for establishing our next result.

\begin{theorem}\label{teo2.4}
Let $u\in C^0(\I;H^s(\R))$, $s>3/2$, be a solution of \eqref{1.0.1} with initial condition $u(0,x)=u_0(x)$. Assume that:
\begin{enumerate}
    \item\label{teo2.3-1} $n$ is odd,
    \item\label{teo2.3-2} For some $\al\in(1/(n+1),1)$, 
    \bb\label{2.2.9}
    |u_0(x)|\sim o(e^{-x})\quad\text{and}\quad |u_0'(x)|\sim O(e^{-\al x})\quad\text{as}\quad x\nnearrow\infty,
    \ee
    and
    \item\label{teo2.3-3} There exists $t_1\in\I$, $t_1>0$, such that
    \bb\label{2.2.10}
    |u(t_1,x)|\sim o(e^{-x})\quad\text{as}\quad x\nnearrow\infty.
    \ee
Then $u\equiv0$.    
\end{enumerate}
\end{theorem}

\begin{remark}
In case $n$ is odd we could replace $O(e^{-\theta x})$ by $O(e^{-\theta |x|})$ in \eqref{2.2.7}--\eqref{2.2.8}, as well as $o(e^{-x})$ by $o(e^{-|x|})$ in \eqref{2.2.9}--\eqref{2.2.10} and $O(e^{-\al |x|})$ in place of $O(e^{-\al x})$, and then let $|x|\rightarrow\infty$. This comes from the following observation: if $n$ is odd, then $v(t,x):=-u(t,-x)$ would also satisfy \eqref{1.0.1} and $v(0,x):=v_0(x)=-u_0(-x)$. Therefore, we can extend the result in theorems \ref{teo2.3} and \ref{teo2.4} for negative values of $x$ through this transformation by applying it to $v$.
\end{remark}

Our next result regards unique continuation of the solutions of \eqref{2.1.2}.

\begin{theorem}\label{teo2.5}
Let $u$ be a solution of \eqref{1.0.1}. Assume that $n$ is odd and there exists a non-empty open rectangle ${\cal R}=(t_0,t_1)\times(x_0,x_1)\subset[0,T)\times\R$ such that $u\big|_{{\cal R}}\equiv0$. 
Then $u\equiv0$.
\end{theorem}

Observe that Theorem \ref{teo2.4} and Theorem \ref{teo2.5} imply local and global properties of the solution in time: while the latter implies that equation \eqref{1.0.1} (for $n$ odd) cannot have compactly supported solutions in $[0,T)\times\R$, the former has a consequence the non-existence of compactly supported solutions for each $t\in(0,T)$, as stated below.

\begin{corollary}\label{cor2.2}
If $n$ is odd and $u$ is a solution of \eqref{1.0.1}, then it cannot be compactly supported in $[0,T)\times\R$.
\end{corollary}

Note that if $n$ is odd and $u$ is a compactly supported solution of \eqref{1.0.1} for some $t>0$, then it satisfies the conditions in Theorem \ref{teo2.4}. Therefore, we have the following straight forward consequence:

\begin{corollary}\label{cor2.3}
If $n$ is odd and $u\in C^0([0,T);H^s(\R))$, $s>3/2$, is a compactly supported solution of \eqref{1.0.1}, then it necessarily is $0$. Conversely, if $u$ does not vanish, then it cannot be compactly supported. In particular, a compactly supported data gives rise to a solution that cannot be compactly supported at any other time. 
\end{corollary}

The next theorem is another unique continuation result for the solutions of \eqref{1.0.1}. 

\begin{theorem}\label{teo2.6}
Assume that $n$ is odd and $u$ is a solution of \eqref{1.0.1}. Suppose that for some $t^\ast\in(0,T)$ we can find distinct points $a$ and $b$, with $a<b$, such that $u_t(t^\ast,a)=u_t(t^\ast,b)=0$ and $u(t^\ast,x)=0$, $x\in[a,b]$. Then $u\equiv0$.
\end{theorem}


Theorems \ref{teo2.1}--\ref{teo2.6} are proved in Section \ref{sec4} whereas in Section \ref{sec3} we present technical results that we need in the demonstration of these theorems. Next, we discuss our achievements in Section \ref{sec5} and present our conclusions in Section \ref{sec6}.

\section{Preliminaries and technical results}\label{sec3}

In this section we present a couple of essential steps needed for the proof of the our main results.

Given a solution $u$ of \eqref{1.0.1}, for each $t\in[0,T)$ we define
\bb\label{3.0.1}
f_t(x):=\f{n}{2}u^{n-1}u_x^2+\f{n(n+3)}{2(n+1)}u^{n+1},
\ee
and
\bb\label{3.0.2}
F_t(x):=\p_x\Lambda^{-2}\ast f_t(x).
\ee

Let us denote the families
\bb\label{3.0.3}
{\cal F}_1=(f_t)_{t\in[0,T)},\quad {\cal F}_2=(F_t)_{t\in[0,T)} 
\ee
and their corresponding sub-families
\bb\label{3.0.4}
{\cal F}^\ast_1=(f_t)_{t\in(0,T)},\quad {\cal F}_2^\ast=(F_t)_{t\in(0,T)}. 
\ee

Note that ${\cal F}_1\subset C^0(\R)$ whereas ${\cal F}_2\subset C^1(\R)$.

\begin{proposition}\label{prop3.1}
Let $u$ be a solution of \eqref{1.0.1} and $f_t(\cdot)$ be the function given by \eqref{3.0.1}. Assume that $n$ is odd. Then we have the equivalences $f_t(\cdot)=0\Leftrightarrow\Lambda^{-2}f_t(x)=0\Leftrightarrow{\cal H}(t)=0$, where ${\cal H}(\cdot)$ is the conserved quantity \eqref{2.2.2}.
\end{proposition}

\begin{proof}
The proof is straightforward once we note that the function \eqref{2.1.1} is positive and $\Lambda^{-2}f_t(x)=(g\ast f_t)(x)$. Therefore, $\Lambda^{-2}f_t(x)=0$ if and only if $f_t(x)\equiv0$. Moreover, we also have $u\equiv0$ if and only if ${\cal H}(\cdot)=\|u(t,\cdot)\|_{H^1(\R)}^2=0$. By \eqref{3.0.1} we realise that if $n$ is odd, then $f_t$ vanishes only when $u$ vanishes.
\end{proof}

\begin{proposition}\label{prop3.2} 
Let $u$ be a solution of \eqref{1.0.1} and assume that $n$ is odd. If there exist numbers $t^\ast$, $x_0$ and $x_1$ such that $\{t^\ast\}\times[x_0,x_1]\subset(0,T)\times\R$, $f_{t^\ast}\big|_{(x_0,x_1)}\equiv0$, $F_{t^\ast}(x_0)=F_{t^\ast}(x_1)$, then ${\cal F}_1={\cal F}_2=\{0\}$. In particular, $u\equiv0$.
\end{proposition}

\begin{proof}
Since $\p_x^2\Lambda^{-2}=\Lambda^{-2}-1$, we have
$$
F_{t^\ast}'(x)=\p_x F_{t^\ast}(x)=(\p_x^2\Lambda^{-2})f_{t^\ast}(x)=\Lambda^{-2}f_{t^\ast}(x)-f_{t^\ast}(x).
$$

Therefore, as long as $x\in(x_0,x_1)$, we conclude that $F_{t^\ast}'(x)=\Lambda^{-2}f_{t^\ast}(x)$. The Fundamental Theorem of Calculus jointly the condition $F_{t^\ast}(x_0)=F_{t^\ast}(x_1)$ yield
$$
0=F_{t^\ast}(x_1)-F_{t^\ast}(x_0)=\int_{x_0}^{x_1}\Lambda^{-2}f_{t^\ast}(x)dx,
$$
which forces $\Lambda^{-2}f_{t^\ast}(x)=0$. The result follows from Proposition \ref{prop3.1}.
\end{proof}

Note that Proposition \ref{prop3.2} has a very strong and beautiful meaning: if we can find a member of ${\cal F}_1^\ast$ that vanishes identically, then ${\cal H}(t)=0$ (see \eqref{2.2.2}) for all $t$, in view of the conservation of the $H^1(\R)-$norm. Moreover, this implies that ${\cal F}_1=\{0\}$, as well as ${\cal F}_2=\{0\}$. On the contrary, if it is possible to find a member of ${\cal F}_1^\ast$ that does not vanish, then any other member of this family will not be identically zero.

\section{Proof of the main results}\label{sec4}

We now combine the results of the previous section to present the demonstration of the theorems announced in Section \ref{sec2}.

{\bf Proof of Theorem \ref{teo2.1}}
A straightforward calculation shows the following formal identities for any $v=v(t,x)$:
\bb\label{v1}
\ba{lcl}
\ds{\p_t(v-v_{xx})+\p_x\left[\f{n+2}{2}v^{n+1}-\f{n}{2}u^{n-1}v_x^2-v^nv_{xx}\right]}&=&\ds{v_t-v_{txx}+\f{(n+1)(n+2)}{2}v^nv_x}\\
\\
\ds{-\f{n(n-1)}{2}v^{n-2}v_x^3+2n v^{n-1}v_xv_{xx}+v^n v_{xxx}}
\ea
\ee
and
\bb\label{v2}
\ba{lcl}
\ds{\p_t\left(\f{v^2+v_x^2}{2}\right)+\p_x\left(\f{v^{n+2}}{n+2}+\f{1-n}{2}v^nv_x^2-v^{n+1}v_{xx}-vv_{tx}\right)=}\\
\\
\ds{v\left(v_t-v_{txx}+\f{(n+1)(n+2)}{2}v^nv_x-\f{n(n-1)}{2}v^{n-2}v_x^3+2n v^{n-1}v_xv_{xx}+v^n v_{xxx}\right)}.
\ea
\ee

Therefore, if $u$ is a solution of \eqref{1.0.1}, the right hand sides of \eqref{v1} and \eqref{v2} vanish. This implies that the left hand sides of these expressions are divergences vanishing on the solutions of \eqref{1.0.1} and, as a consequence, they are conserved currents for the equation. Moreover, $Q=1$ and $Q=u$ are the corresponding characteristics.

To the uniqueness of these conserved currents, let $\phi=\phi(t,x,u)$ and
$$
E:=u_t-u_{txx}+\f{(n+1)(n+2)}{2}u^nu_x-\f{n(n-1)}{2}u^{n-2}u_x^3+2n u^{n-1}u_xu_{xx}+u^n u_{xxx}.
$$

By \cite[Theorem 4.7]{olverbook}, $Q E$ is a divergence if and only if $E_u(Q E)=0$, where $E_u$ is the Euler-Lagrange operator (see \cite[Def. 4.3]{olverbook}). Therefore, solving the equation $E_u(Q E)=0$ we find $Q=c_1+c_2u$, where $c_1$ and $c_2$ are two arbitrary constants, meaning that the unique characteristics are just $Q=1$ and $Q=u$.
\hfill$\square$

{\bf Proof of Corollary \ref{cor2.1}}
Assuming \eqref{2.2.1}, if $u$ is either non-negative or non-positive, then $u=\pm|u|$. Therefore, from \eqref{2.2.1} we conclude that ${\cal H}_1(t)=\pm\|u(t,\cdot)\|_1$.

Let us prove that \eqref{2.2.1} and \eqref{2.2.2} are constant. If we choose $v=u$ and integrate \eqref{v1} and \eqref{v2} over $\R$, we obtain, respectively,
$$
\f{d}{dt}\int_\R(u-u_{xx})dx=\f{d}{dt}\int_\R u dx=0
$$
and
$$
\f{d}{dt}\int_\R\f{u(t,x)^2+u_x(t,x)^2}{2}dx=0.
$$
The relations above imply the desired results. \hfill$\square$

We now prove theorems \ref{teo2.5} and \ref{teo2.6} and then the remaining results. Note that the conditions in theorems \ref{teo2.5} and \ref{teo2.6} imply that $f_t(\cdot)$ is non-negative.

{\bf Proof of Theorem \ref{teo2.5}.} From \eqref{3.0.2} and \eqref{2.1.2} we have
\bb\label{4.0.1}
F_t(x)=-(u_t+u^nu_x)(t,x),\quad t>0.
\ee

Therefore, if $u$ vanishes on ${\cal R}$, then $F_t(x)=0$ provided that $(t,x)\in{\cal R}$. Let us choose $t^\ast\in(t_0,t_1)$ and real numbers $a$ and $b$ such that $x_0<a<b<x_1$. Then
$$
F_{t^\ast}(b)=F_{t^\ast}(a)=0.
$$
Since $f_{t^\ast}(x)\geq0$ and $f_{t^\ast}\big|_{[a,b]}\equiv0$, then Proposition \ref{prop3.2} implies the result.\hfill$\square$

{\bf Proof of Theorem \ref{teo2.6}.} Let us assume that we could find $t^\ast\in(0,T)$ and an interval $[a,b]$ such that $u(t^\ast,x)=0$, $x\in[a,b]$, and $u_t(t^\ast,a)=u_t(t^\ast,b)=0$. By \eqref{4.0.1} we conclude that $F_{t^\ast}(a)=F_{t^\ast}(b)=0$. Also, we have $f_{t^\ast}(x)=0$, $x\in[a,b]$. The result is again a consequence of Proposition \ref{prop3.2}.
\hfill$\square$

Before proving theorems \ref{teo2.3} and \ref{teo2.4}, let us define
\bb\label{4.0.2}
M:=\sup_{t\in\I}\|u(t,\cdot)\|_{H^s(\R)}+\sup_{t\in\I}\|u_x(t,\cdot)\|_{H^s(\R)}<\infty.
\ee

Also, in what follows $p$ denotes a positive integer, and we would like to note that $\p_xg\in L^1(\R)$, where $g$ is given by \eqref{2.1.1}, and $u^{n-1}u_x^2,\,u^{n+1}\in L^1(\R)\cap L^\infty(\R)$ in view of the Sobolev Embedding Theorem \cite[Proposition 1.3]{taylor} and the Algebra Property \cite[page 320, Exercise 6]{taylor}. Therefore, we have
$$
{\cal F}_1\subset L^1(\R)\cap L^\infty(\R)\cap C^0(\R)\quad\text{and}\quad {\cal F}_2\subset L^1(\R)\cap L^\infty(\R)\cap C^1(\R)
$$
where the families ${\cal F}_1$ and ${\cal F}_2$ are given by \eqref{3.0.1}.

In addition, noticing that if $h\in L^1(\R)\cap L^\infty(\R)$, then $h\in L^r(\R)$, for all $1\leq r\leq \infty$, from which we have
$$\lim_{r\rightarrow\infty}\|h\|_r=\|h\|_\infty.$$

Finally, in the next pages we make use of two claims, whose demonstrations are left to the Appendix. This is done in order to allow a more fluid reading and understanding of the demonstration and avoid moderately lengthy technical parts of the theorems.

{\bf Proof of Theorem \ref{teo2.3}} Our strategy for proving Theorem \ref{teo2.3} is the following: we construct a sequence of piecewise smooth functions $\phi_N(x)$, such that $\phi_N(x)\rightarrow e^{\theta x}$ pointwisely for $x>0$, $\|\phi_N(\cdot)u(t,\cdot)\|_\infty$ and $\|\phi_N(\cdot)u_x(t,\cdot)\|_\infty$ are bounded, for each $N$. Similarly as in \cite{him-cmp}, we show the existence of $L\geq0$ such that $|\phi_N(x)u(t,x)|+|\phi_N(x)u_x(t,x)|$ is bounded by $L$. Then, allowing $N\rightarrow\infty$ and using \eqref{2.2.7} we conclude \eqref{2.2.8}.

Let $N$ be a positive integer and consider the family $\phi_N:\R\rightarrow\R$, given by
$$
\phi_N(x)=
\left\{
\ba{lcl}
1,&\text{if}&x\leq0,\\
\\
e^{\theta x},&\text{if}&0<x<N,\\
\\
e^{\theta N},&\text{if}&x\geq N.
\ea
\right.
$$

Clearly $(\phi_N)_{N\in\mathbb{N}}\subset C^0(\R)$ and, for each $N$, $\phi_N$ is piecewise smooth and $\phi_N'(x)\leq\phi_N(x)$ almost everywhere.

Multiplying \eqref{2.1.2} by $\phi_N(\phi_Nu)^{2p-1}$ and integrating with respect to $x$ over $\R$, we obtain
$I_1=-I_2-I_3$, where
$$
I_1=\|\phi_N(\cdot) u(t,\cdot)\|_{2p}^{2p-1}\f{d}{dt}\|\phi_N(\cdot) u(t,\cdot)\|_{2p},
$$
$$
I_2=\int_\R(\phi_N(x) u(t,x))^{2p}u(t,x)^{n-1}u_x(t,x)dx\quad\text{and}\quad I_3=\int_\R(\phi_N(x) u(t,x))^{2p-1}\phi_N (x)F_t(x)dx.
$$

{\bf Claim 1:} $|I_2|\leq M^n\|\phi_N(\cdot)u(t,\cdot)\|_{2p}^{2p}$ and $|I_3|\leq \|\phi_N(\cdot)u(t,\cdot)\|_{2p}^{2p-1}\|\phi_N(\cdot) F_t(\cdot)\|_{2p}$.

Since $I_1=-I_2-I_3$, the estimates in the claim gives the differential inequality
$$
\f{d}{dt}\|\phi_N(\cdot)u(t,\cdot)\|_{2p}\leq M^n\|\phi_N(\cdot)u(t,\cdot)\|_{2p}+\|\phi_N(\cdot) F_t(\cdot)\|_{2p},
$$
that jointly with Gr\"onwall's inequality and letting $p\rightarrow\infty$, yield
\bb\label{4.0.5}
\|\phi_N(\cdot)u(t,\cdot)\|_{\infty}\leq\left(\|\phi_n(\cdot)u_0(\cdot)\|_\infty+\int_0^t e^{-M^nt}\|\phi_N(\cdot) F_\tau(\cdot)\|_\infty d\tau\right)e^{M^nt}.
\ee

Differentiating \eqref{2.1.2} with respect to $x$, we have
\bb\label{neweq}
u_{tx}+nu^{n-1}u_x^{2}+u^nu_{xx}+\p_x F_t(x)=0.
\ee

Now we multiply \eqref{neweq} by $\phi_N(\phi_Nu_x)^{2p-1}$ to obtain the inequality $I_4\leq |I_5|+|I_6|+|I_7|$, where
$$
\ba{lcl}
I_4&=&\ds{\|\phi_N (\cdot)u_x(t,\cdot)\|_{2p}^{2p-1}\f{d}{dt}\|\phi_N (\cdot)u_x(t,\cdot)\|_{2p},}\\
\\
I_5&\leq&\ds{\int_\R n u(t,x)^{n-1}(\phi_N(x)u_x(t,x))^{2p}u_x(t,x)dx},\\
\\
I_6&\leq& \ds{\int_\R u(t,x)^n(\phi_N(x)u_x(t,x))^{2p-1}(\phi_N(x)u_{xx}(t,x))dx},\\
\\
I_7&\leq&\ds{\int_\R \phi_N(x)(\phi_N(x)u_x(t,x))^{2p-1}\p_x F_t(x)dx},
\ea
$$
and we have our last claim.

{\bf Claim 2:} $|I_5|\leq n M^n \|\phi_N (\cdot)u_x(t,\cdot)\|_{2p}^{2p}$, $|I_6|\leq (n+1)M^n \|\phi_N (\cdot)u_x(t,\cdot)\|_{2p}^{2p}$ and $|I_7|\leq\|\phi_N (\cdot)u_x(t,\cdot)\|_{2p}^{2p-1}\|\phi_N(\cdot) \p_x F_t(\cdot)\|_{2p}$.

Therefore, we have the differential inequality
$$
\f{d}{dt}\|\phi_N (\cdot)u_x(t,\cdot)\|_{2p}\leq 2(n+1)M^n\|\phi_N (\cdot)u_x(t,\cdot)\|_{2p}+\|\phi_N(\cdot)\p_x F_t(\cdot)\|_{2p},
$$
which yields, after using the Gr\"onwall inequality, and taking $p\rightarrow\infty$ again, we have
\bb\label{4.0.6}
\|\phi_N (\cdot)u_x(t,\cdot)\|_{\infty}\leq\left(\|\phi_N(\cdot) u_0'(\cdot)\|_\infty+\int_0^te^{-2(n+1)M\tau}\|\phi_N(\cdot) \p_xF_\tau(\cdot)\|_\infty d\tau
\right)e^{2(n+1)Mt}.
\ee

We not recall that there exists a positive constant $c$ (see Appendix \ref{apC}) such that
\bb\label{4.0.7}
\ba{lcl}
\ds{\|\phi_N(\cdot) F_\tau(\cdot)\|_\infty}&\leq&\ds{c\left(\|\phi_N (\cdot)u(t,\cdot)\|_\infty+\|\phi_N (\cdot)u_x(t,\cdot)\|_\infty\right)}\\
\\
\ds{\|\phi_N(\cdot) \p_xF_\tau(\cdot)\|_\infty}&\leq&\ds{c\left(\|\phi_N (\cdot)u(t,\cdot)\|_\infty+\|\phi_N (\cdot)u_x(t,\cdot)\|_\infty\right).}
\ea
\ee

Also, we recall that $e^{-M^nt}, e^{-2(n+1)t}\in(0,1)$. Adding \eqref{4.0.5} to \eqref{4.0.6}, and taking into account \eqref{4.0.7} and that $t\in\I=[0,T_0]$, we have
$$
\ba{lcl}
\ds{\|\phi_N (\cdot)u(t,\cdot)\|_\infty+\|\phi_N (\cdot)u_x(t,\cdot)\|_\infty}&\leq&\ds{ K\Big[\left(\|\phi_N(\cdot) u_0(\cdot)\|_\infty+\|\phi_N(\cdot) u'_0(\cdot)\|_\infty\right)}\\
\\
&&\ds{+\int_0^t\left(\|\phi_N (\cdot)u(t,\cdot)\|_\infty+\|\phi_N (\cdot)u_x(t,\cdot)\|_\infty\right)d\tau\Big]},
\ea
$$
for some $K=K(M,n,t_0)>0$. which, after application of Gr\"onwall inequality, yields
\bb\label{4.0.8}
\|\phi_N (\cdot)u(t,\cdot)\|_\infty+\|\phi_N (\cdot)u_x(t,\cdot)\|_\infty\leq\kappa\left(\|\phi_N(\cdot) u_0(\cdot)\|_\infty+\|\phi_N(\cdot) u'_0(\cdot)\|_\infty\right),
\ee
where $\kappa:=KT_0$.

On the other hand, we have
$$
|\phi_N(x)u_0(x)|\leq |\max\{1,e^{\theta x}\}u_0(x)|\leq\|\max\{1,e^{\theta \cdot}\}u_0(\cdot)\|_\infty=:L_1/2
$$
and
$$
|\phi_N(x)u_0'(x)|\leq |\max\{1,e^{\theta x}\}u'_0(x)|\leq\|\max\{1,e^{\theta \cdot}\}u'_0(\cdot)\|_\infty=:L_2/2
$$

Let $L:=\max\{\kappa L_1,\kappa L_2\}$. The inequalities above, jointly with \eqref{4.0.8}, implies that
\bb\label{4.0.9}
\|\phi_N (\cdot)u(t,\cdot)\|_\infty+\|\phi_N (\cdot)u_x(t,\cdot)\|_\infty\leq L.
\ee

Note that $L$ does not depend on $N$ and if $x>0$, taking the limit $N\rightarrow\infty$, we conclude that 
$$\phi_N(x)u(t,x)\rightarrow e^{\theta x}u(t,x)\quad\text{and}\quad\phi_N(x)u_x(t,x)\rightarrow e^{\theta x}u_x(t,x).$$

Therefore, if $N\rightarrow\infty$ and $x>0$, \eqref{4.0.9} gives
$$
|e^{\theta x}u(t,x)|+|e^{\theta x}u_x(t,x)|\leq L,
$$
which implies the result. In particular, note that if $u_0(x_0)\neq0$ for some $x_0$, then \eqref{4.0.9} implies that $L$ is positive.
\hfill$\square$

{\bf Proof of Theorem \ref{teo2.4}:} Integrating \eqref{4.0.1} with respect to $t$ we have
$$
-\int_0^{t_1} F_\tau(x)d\tau=u(t_1,x)-u_0(x)+\int_0^{t_1}(uu^n)(\tau,x)d\tau=:J.
$$

By Theorem \ref{teo2.3}, we conclude that $(uu^n)(\tau,x)\sim O(e^{-\al(n+1) x})$. Moreover, since $1/(n+1)<\al<1$, then $1-\al(n+1)<0$. Thus,
$$
\lim_{x\rightarrow\infty}\f{|(uu^n)(t,x)|}{e^{-x}}=\lim_{x\rightarrow\infty}\Big(\underbrace{\f{|(uu)   ^n(t,x)|}{e^{-\al(n+1)x}}}_{bounded}\,\,\underbrace{\f{e^{-\al(n+1)x}}{e^{-x}}}_{\rightarrow0}\Big)=0,
$$
meaning that $(uu^n)(\tau,x)\sim o(e^{-x})$. Therefore, it is immediate that $J\sim o(e^{-x})$.

We will prove the result by contradiction: we assume that $u$ is not identically zero. Then we show that if this is true, we can find a positive constant $K$ such that 
$$-\int_0^{t_1} F_\tau(x)d\tau\geq Ke^{-x}$$
which contradicts $J\sim o(e^{-x})$.

We begin by recalling that $F_t(x)=(\p_xg)\ast f_t(x)$, where $f_t(x)$ is given by (see \eqref{3.0.1})
$$
f_t(x)=\f{n}{2}u(t,x)^{n-1}\left(u_x(t,x)^2+\f{n+3}{n+1}u(t,x)^2\right)\sim o(e^{-(n+1)\al x}).
$$

Let us define 
$$\rho(x):=\int_0^{t_1}f_\tau(x)d\tau.$$

Then we have the following identity:
\bb\label{4.0.10}
-\int_0^{t_1}F_\tau(x)d\tau=(\p_xg)\ast\rho(x).
\ee

We also observe that
\bb\label{4.0.11}
(\p_xg )\ast\rho(x)=-\f{1}{2}e^{-x}\int_{-\infty}^x e^y\rho(y)dy+\f{1}{2}e^{x}\int_{x}^\infty e^{-y}\rho(y)dy,
\ee
and
$$
e^x\int_x^\infty\rho(y)dy=o(1)e^x\int_x^\infty e^{-(n+2)\al y}dy\sim o(1)e^{-(n+1)x}=o(e^{-(n+1)x})\sim o(e^{-x}).
$$

On the other hand, if $n$ is odd, then $\rho(y)>0$ (recall that $u\not\equiv0$) and 
$$x\mapsto \int_{-\infty}^x e^y\rho(y)dy$$
is an increasing function, so that for $x$ sufficiently large, we can find a constant $k>0$ such that
$$
 \int_{-\infty}^x e^y\rho(y)dy\geq k>0.
$$

The latter result, jointly with \eqref{4.0.11} and \eqref{4.0.10} imply that
$$
-\int_0^{t_1}F_\tau(x)d\tau\geq \f{k}{2} e^{-x},
$$
which shows the desired contradiction.
\hfill$\square$

\section{Discussion}\label{sec5}

To the best of the author's knowledge, the Lie point symmetries of he equation \eqref{1.0.1} has not been previously investigated, as well as a classification of its lower order conservation laws. The Lie algebra of the Lie symmetries of \eqref{1.0.1} is just the same as the CH equation  (compare \cite[Proposition 1.1]{anco} with Theorem \ref{teo2.2}), which is somewhat expected in view of the structure of the equation. In regard to conservation laws, we showed that the zeroth order characteristics are $\phi_1=1$ and $\phi_2=u$, that correspond to conserved currents of zero and first orders, respectively. We note that these characteristics are just the same admitted by the CH equation, see \cite[Proposition 2.1]{anco} and \cite[Theorem 2.1]{anco}, respectively.

We observe that from the symmetries we found some solutions of \eqref{1.0.1}. More precisely, from the scaling $(t,x,u)\mapsto(e^{n \epsilon}t,x,e^{-\epsilon}u)$ we obtained two solutions for \eqref{1.0.1}: $u(t,x)=e^{x}/t^{1/n}$ and $u(t,x)=e^{-x}/t^{1/n}$. These solutions, of course, do not conserve both \eqref{2.2.1} and \eqref{2.2.2}. Another sort of solution coming from the Lie symmetries are the traveling waves. In particular, these solutions satisfy the relation \eqref{2.2.5}. Although \eqref{2.2.5} was obtained, in principle, for classical solutions, the fact that it may have a singularity lead us naturally to consider solutions in the week sense and, in particular, we show that the equation has peakon solutions for each $n$. Moreover, as shown in \cite{hak-ijns}, these peakons are orbitally stable solutions for \eqref{1.0.1} for each positive integer $n$.

As far as the author knows, the first (and unique) work considering asymptotic properties of \eqref{1.0.1} with $n\neq1$ is \cite{zhou-jmaa}, in which the ideas introduced by Brandolese \cite{bran-imrn} for the studying the CH equation for \eqref{1.0.1} were extended. In the present work we study these properties in a different way, using the ideas from \cite{him-cmp} for proving theorems \ref{teo2.3} and \ref{teo2.4}.

We note that our results are mostly concerned when $n$ is odd. Such condition comes from the necessity that the function \eqref{3.0.1} be non-negative. Essentially, this also explains why the assumption that $n$ is odd appears in \cite[Theorem 1.3]{zhou-jmaa}. We could replace this hypothesis by assuming that $u$ is either non-negative or non-positive, however, whether \eqref{1.0.1} has non-negative or non-positive solutions for $n>1$ seems to be an open question. For the CH equation this is proved by using a diffeomorphism constructed from its solutions \cite[Theorem 3.1]{const2000-1}. Although we could find a similar diffeomorphism for \eqref{1.0.1}, we do not have a similar result as \cite[Lemma 3.2]{const2000-1}. Moreover, the apparent impossibility of the extension\footnote{Or the fact that so far we (or perhaps, the author,) do not know how to reach to a similar result.} of such result to \eqref{1.0.1} make us also unable to proceed similarly as in \cite[Theorem 1.4]{him-cmp} to prove that \eqref{1.0.1} does not have compactly supported solutions corresponding to an initial data compactly supported.

Finally, we also reported some unique continuation results for \eqref{1.0.1} based on recent ideas \cite{linares}, see also \cite{freire-jpa,freire-cor}. As pointed out in \cite{raspa-mo,freire-jpa,freire-arxiv}, these ideas are essentially geometric, see \cite{raspa-mo,freire-jpa,freire-arxiv} for a better discussion.

\section{Conclusion}\label{sec6}

In this paper we studied persistence and asymptotic properties of an extension of the CH equation, given by \eqref{1.0.1}. Our main results are theorems \ref{teo2.3}--\ref{teo2.6}. In particular, we complement the results in \cite{zhou-jmaa} by showing that \eqref{1.0.1} does not have compactly supported solutions, as well as we extend to \eqref{1.0.1} (for odd $n$) the results proved in \cite{linares}, see also \cite{freire-jpa,freire-arxiv}, regarding unique continuation of solutions.

\section*{Acknowledgements}

I. L. Freire is grateful to FAPESP for financial support (grant nº 2020/02055-0).

\section*{Data availability}

Data sharing is not applicable to this article as no new data were created or analyzed in this study.

\appendix











\section{Proof of Claim 1}

Note that
$$
|I_2|\leq\|u(t,\cdot)\|_\infty^{n-1}\|u_x(t,\cdot)\|_\infty\int_\R(\phi_N(x)u(t,x))^{2p}dx\leq\|u(t,\cdot)\|_\infty^{n-1}\|u_x(t,\cdot)\|_\infty\|\phi_N(\cdot)u(t,\cdot)\|_{2p}^{2p}
$$
and, from Hölder's inequality, we have
$$
|I_3|\leq\int_\R|\phi_N(x)u(t,x)|^{2p-1}|\phi_N(x)F_t(x)|dx\leq\|\phi_N(\cdot)u(t,\cdot)\|_{2p}^{2p-1}\|\phi_N(\cdot)F_t(\cdot)\|_{2p}.
$$

The result is again a consequence of the inequalities above and \eqref{4.0.2}.

\section{Proof of Claim 2}

It is straightforward that
$$
|I_5|\leq n\|u(t,\cdot)\|_{\infty}^{n-1}\|u_x(t,\cdot)\|_\infty\int_\R (\phi_N(x)u_x(t,x))^{2p}dx
$$
and
$$\int_\R (\phi_N(x)u_x(t,x))^{2p}dx\leq \|\phi_N(\cdot)u_x(t,\cdot)\|_{2p}^{2p}.$$

In regard to $I_6$, integration by parts gives
$$
I_6=-\f{n}{2p}\int_\R u(t,x)^{n-1}u_{x}(t,x)(\phi_N(x)u_x(t,x))^{2p}dx-\int_\R u(t,x)^n(\phi_N(x)u_x(t,x))^{2p-1}\phi'_n(x)u_x(t,x)dx,
$$
which implies after noticing that $\phi_N(x)\leq\phi_N(x)$ a.e., 
$$
|I_6|\leq (n+1)\Big(\|u(t,\cdot)\|_\infty^{n-1}\|u_x(t,\cdot)\|_\infty+\|u(t,\cdot)\|_\infty^n\Big)\|\phi_N(\cdot)u_x(t,\cdot)\|_{2p}^{2p}.
$$

With respect to $I_7$, using the Hölder inequality we have
$$
|I_7|\leq \int_\R|\phi_N(x)u_x(t,x)^{2p-1}||\phi_N(x)\p_x F_t(x)|dx\leq \|\phi_N(\cdot)u_x(t,\cdot)\|_{2p}^{2p-1}\|\phi_N(\cdot)\p_x F_t(\cdot)\|_{2p}.
$$

The inequalities above jointly with \eqref{4.0.2} complete the demonstration of Claim 2.

\section{Proof of estimates \eqref{4.0.7}}\label{apC}

We begin by noticing that $|\phi_N(x)F_\tau(x)|=|\phi_N(x)((\p_xg)\ast f_\tau)(x)|$. By \eqref{3.0.1} and \eqref{4.0.2}, and \cite[Eq. (2.31)]{him-cmp}, which shows that
$$
\phi_N(x)\int_\R \f{e^{-|x-y|}}{\phi_N(y)}dy=\f{4}{1-\theta}=:c_0,
$$
we have
$$
\ba{lcl}
\ds{|\phi_N(x)F_\tau(x)|}&\leq&\ds{k\|u\|_\infty^{n-1}\phi_N(x)\int_\R e^{-|x-y|}\left(u^2+u_x^2\right)(\tau,y)dy}\\
\\
&\leq&\ds{kM^{n-1}\Bigg[\left(\phi_N(x)\int_\R \f{e^{-|x-y|}}{\phi_N(y)}dy\right)\|\phi_N(\cdot)u_x(t,\cdot)\|_\infty\|u_x(t,\cdot)\|_\infty}\\
\\
&&\ds{\left(\phi_N(x)\int_\R \f{e^{-|x-y|}}{\phi_N(y)}dy\right)\|\phi_N(\cdot)u(t,\cdot)\|_\infty\|u(t,\cdot)\|_\infty\Bigg]}\\
\\
&\leq&\ds{ kM^nc_0(\|\phi_N(\cdot)u(t,\cdot)\|_\infty+\|\phi_N(\cdot)u_x(t,\cdot)\|_\infty).}
\ea
$$

Since $\p_x^2g=g-\delta$, where $\delta$ is the Dirac $\delta$ distribution, the same procedure shows the second estimate in \eqref{4.0.7} and for this reason its proof is omitted.
\end{document}